\documentclass[11pt,reqno]{amsart}
\usepackage{amsmath,amssymb,amsthm,mathtools,mathrsfs}
\usepackage[margin=1in]{geometry}
\usepackage[hidelinks]{hyperref}
\hypersetup{pdftitle={Limiting particle density and vacuum concentration in the one-dimensional compressible Navier-Stokes-Vlasov system},pdfauthor={Hanyu Teng}}
\usepackage{microtype,needspace}
\newtheorem{theorem}{Theorem}[section]
\newtheorem{lemma}[theorem]{Lemma}
\newtheorem{proposition}[theorem]{Proposition}
\newtheorem{corollary}[theorem]{Corollary}
\theoremstyle{remark}
\newtheorem{remark}[theorem]{Remark}
\numberwithin{equation}{section}
\newcommand{\R}{\mathbb R}
\newcommand{\T}{\mathbb T}
\newcommand{\dd}{\,\mathrm d}
\newcommand{\I}{\mathcal I}
\newcommand{\E}{\mathcal E}
\newcommand{\av}[1]{\langle #1\rangle}
\DeclareMathOperator{\TV}{TV}
\title[Limiting particle density and vacuum concentration]
{Limiting particle density and vacuum concentration in the one-dimensional compressible Navier--Stokes--Vlasov system}

\author[H. Teng]{Hanyu Teng}
\address{School of Mathematical Sciences, Xiamen University, Xiamen, 361005, P.R. China.}
\email{19020260158071@stu.xmu.edu.cn}
\date{November 26, 2026}
\keywords{Compressible Navier--Stokes--Vlasov system, vacuum, finite-energy weak solutions, specular reflection, asymptotic profile, renormalization}
\subjclass[2020]{35Q30,35Q83,35B40,76N10}
\begin{document}
\begin{abstract}
We construct global finite-energy weak solutions to the one-dimensional compressible Navier--Stokes--Vlasov system on a reflecting interval, allowing vacuum and unbounded initial fluid density. The solutions satisfy the energy equality at every time. At positive mean fluid density $r$, a conserved quantity along particle paths determines their limiting positions from the initial data. For bounded particle data of finite mass and kinetic energy, the spatial density converges strongly in $L^3$, and the exponent is sharp. If the initial vacuum has positive length, its length times $e^{r^\gamma t/\nu}$ has a positive essential limit. A weighted identity in the fluid mass coordinate determines the normalized limiting vacuum measure from the data.
\end{abstract}
\maketitle

\section{Introduction}\label{intro}
We consider a compressible fluid coupled to a particle distribution on an interval $I=(0,L)$. The fluid density $\rho$, the fluid velocity $u$ and the particle density in phase space $f=f(t,x,v)$ satisfy
\begin{equation}\label{sys}
\begin{cases}
 \rho_t+(\rho u)_x=0,\\
 (\rho u)_t+(\rho u^2+\rho^\gamma)_x=\nu u_{xx}-\rho\int_\R(u-v)f\dd v,\\
 f_t+vf_x+\partial_v\big(\rho(u-v)f\big)=0.
\end{cases}
\end{equation}
Here $\gamma>1$ and $\nu>0$. The fluid velocity vanishes at each endpoint, where the particles reflect from $v$ to $-v$. We write $m=\rho u$, $n(t,x)=\int_\R f(t,x,v)\dd v$, and assume
\begin{equation}\label{data}
 0\le\rho_0\in L^\gamma(I),\quad \int_I\frac{|m_0|^2}{\rho_0}\dd x<\infty,
 \qquad 0\le f_0\in L^\infty(I\times\R),\quad
 \iint_{I\times\R}(1+v^2)f_0\dd x\dd v<\infty .
\end{equation}
The kinetic energy density is zero where $\rho_0=m_0=0$ and infinite where $\rho_0=0$, $m_0\ne0$.

For the same drag on the torus, Li and Shou \cite[Theorems 1.1--1.2]{LS} proved exponential alignment for initial fluid density in $W^{1,\infty}$ bounded away from zero, initial velocity in $H^1$ and bounded particle data of compact velocity support. Their Wasserstein estimate compares $f(t)$ with $n(t,x)\dd x\otimes\delta_{u_s}$, where $u_s$ is the common limiting velocity. It leaves open the limiting spatial density. A similar question arises for the fluid. Convergence of $\rho$ to its mean in $L^2$ gives little information about the location or length of any remaining vacuum.

Both distributions can be described in terms of the initial mass to the left of a point. When the mean fluid density $r$ is positive, set
\begin{equation}\label{masscoord}
 r=\frac{1}{L}\int_0^L\rho_0(x)\dd x,\qquad
 S_0(a)=\frac{1}{r}\int_0^a\rho_0(x)\dd x .
\end{equation}
Write $Z_t=\{x\in I:\rho(t,x)=0\}$ for the vacuum set and $|Z_t|$ for its length.

\Needspace{7\baselineskip}
\begin{theorem}\label{main}
Under \eqref{data}, system \eqref{sys} has a global finite-energy weak solution satisfying the energy equality at every time. Each solution constructed in Appendix~\ref{construction} with $r>0$ has the following properties.

\smallskip\noindent\textup{(i)} The particle density converges strongly in $L^3(I)$ to $n_\infty$. For every $\varphi\in C([0,L])$,
\begin{equation}\label{physprofile}
 \int_0^L\varphi(x)n_\infty(x)\dd x
 =\int_0^L\!\int_\R \varphi\!\left(\operatorname{dist}
       (S_0(a)+v/r,\,2L\mathbb Z)\right) f_0(a,v)\dd v\dd a .
\end{equation}
Here $\operatorname{dist}(y,2L\mathbb Z)$ folds a point of the covering line back into the physical interval.

\smallskip\noindent\textup{(ii)} There is a finite singular measure $\lambda_\infty$ on $[0,L]$ such that, as $t\to\infty$ through a set of full measure,
\[
 e^{r^\gamma t/\nu}\int_{Z_t}\varphi(x)\dd x
      \longrightarrow\int_{[0,L]}\varphi\dd\lambda_\infty
      \qquad\text{for every }\varphi\in C([0,L]).
\]
In particular, writing $V_*=\lambda_\infty([0,L])$, we have
\begin{equation}\label{vaclim}
 \operatorname*{ess\,lim}_{t\to\infty} e^{r^\gamma t/\nu}|Z_t|=V_*,
 \qquad V_*>0\quad\Longleftrightarrow\quad |Z_0|>0 .
\end{equation}
When $|Z_0|>0$, the initial data and $\nu$ determine $\lambda_\infty/V_*$. The fraction in a Borel set $A$ is the weighted initial length of those $a\in Z_0$ for which $S_0(a)\in A$, as given by \eqref{vacnormalized}.
\end{theorem}

The construction applies to every datum in \eqref{data}. We approximate by smooth positive fluid densities as in \eqref{approx}. For $r>0$, every subsequential limit obtained in Appendix~\ref{construction} has the stated particle profile and, when $|Z_0|>0$, the same normalized vacuum measure. Uniqueness of the coupled evolution in the finite-energy class is open. For fixed $\nu,\gamma,L$, the bound $\|\rho(t)\|_\infty\le C(1+t^{-1/\gamma})$ uses only upper bounds on the initial energy and the two masses, and remains uniform as the approximating fluid mass tends to zero. At $r=0$ the fluid is at rest and the particles undergo free specular transport.

For example, on $(0,1)$ take $\rho_0=4/3$ outside $(1/4,1/2)$ and zero inside. The mean density is one, and the entire gap has mass coordinate $1/3$. It gives an atom of $\lambda_\infty$ at $1/3$. A particle in the gap with initial velocity $1$ has limiting position $2/3$ after reflection. More generally, \eqref{physprofile} depends on $\rho_0$ and $f_0$ alone. The relative masses of the vacuum atoms can also depend on the viscosity, as the two-gap example in Section~\ref{vacuum} shows.

\subsection*{Method of proof and related results}
Let $S(t,x)=r^{-1}\int_0^x\rho(t,y)\dd y$. Along the fluid, the continuity equation gives $S_t+uS_x=0$. Along an unfolded particle trajectory $(X,V)$, its derivative is $\rho(V-u)/r$, precisely the opposite of $\dot V/r$. Hence
\begin{equation}\label{twolabels}
 S(t,X(t))+V(t)/r=S_0(X(0))+V(0)/r .
\end{equation}
As $\rho\to r$ and $V\to0$, \eqref{twolabels} gives the limiting position. For strong convergence of the density, we also estimate the derivative of $X$ in the conserved label (Lemma~\ref{jaclemma}) and control the tail of the initial velocity moment. Kim \cite{Kim} used a first-order reduction for a different alignment law, and Jabin \cite{Jabin} derived a limiting spatial measure from kinetic-energy decay. For incompressible Vlasov--Navier--Stokes flows, periodic profiles are studied in \cite{HKMM,DS}, and Han-Kwan \cite{HK} obtained polynomial decay for small data in $\R^3$. On a bounded domain, Ertzbischoff, Han-Kwan and Moussa \cite{EHM} proved velocity concentration for small data with absorbing walls, where part or all of the initial particle mass can escape. Here reflection preserves all the mass, and \eqref{physprofile} locates it from the initial data. Related alignment results treat density-independent drag \cite{Choi}, inhomogeneous incompressible flow \cite{LSZ}, and locally averaged drag \cite{CS}.

For the uncoupled fluid, Tan, Wang and Zhang \cite[Corollary 1.1]{TWZ} located a prescribed vacuum bubble by the fluid mass initially to its left. In the coupled problem the same coordinate locates the vacuum, but it gives only its position. To determine its weight, we use the effective fluid velocity of Li and Shou \cite[(1.18)]{LS}. On a reflecting interval its correction is simply the particle mass to the left of $x$. Integrating the corrected momentum equation expresses $\nu u_x$ as a material derivative, the local pressure and a spatial average. On the vacuum the local pressure vanishes. The remaining average multiplies all vacuum lengths by the same factor, which cancels when their proportions are taken.\footnote{OpenAI's ChatGPT assisted in exploring this derivation and checking the two-gap calculation in Section~\ref{vacuum}.} Li and Shou's companion paper \cite{LSChinese} establishes fluid relaxation for regular solutions admitting vacuum under a priori bounds.

In our earlier work \cite[Theorem 1.4]{Teng}, Besov estimates and boundary cutoffs gave renormalization and the energy equality for the compressible system in two and three dimensions under upper and positive lower density bounds. Here the same two issues arise for data admitting vacuum. The conserved particle quantity reduces the kinetic equation to transport with a spatially Lipschitz velocity. For the fluid, the one-dimensional product estimate in $H^1$ justifies the kinetic-energy chain rule at vacuum. The velocity cutoff may again be chosen even, as required by specular reflection.

Section~\ref{lemmas} gives the transport identities. Sections~\ref{particles} and \ref{vacuum} prove the two parts of Theorem~\ref{main}. Appendix~\ref{construction} constructs the weak solutions and proves the energy equality. Appendix~\ref{general-periodic} gives the limiting distributions in a moving frame for periodic data.

\section{Transport identities}\label{lemmas}
Reflect $I$ to the torus $\T_\ell$ of length $\ell=2L$. Then $\rho$ is even, $u$ is odd, and $f(t,-x,-v)=f(t,x,v)$. On the covering line the cumulative mass and its periodic part are
\[
 M(t,x)=\int_0^x\rho(t,s)\dd s,\qquad P(t,x)=M(t,x)-rx.
\]
Thus $M_x=\rho$, $M_t=-m$ by continuity and parity, and $S=M/r$ if $r>0$. We write $[y]_\ell$ for reduction modulo $\ell$, $\av{h}=\ell^{-1}\int_0^\ell h(x)\dd x$ for the spatial average, and $\I[h]$ for the unique periodic function of mean zero with $(\I[h])_x=h-\av{h}$. In particular $P=\I[\rho]$ by symmetry. Norms without a domain are taken over one period.

\begin{lemma}[A conserved particle label]\label{masslemma}
Let $M\in C([0,T]\times\R)$ with $M_x=\rho\ge0$, $M_t=-m$ in $L^1((0,T)\times K)$ on each bounded interval $K$, $\|\rho(t)\|_\infty\in L^1(0,T)$ and $M(t,x+\ell)=M(t,x)+r\ell$. For a bounded distributional particle density with prescribed initial trace and flux $(m-\rho v)f$ locally integrable up to zero, set $w=v+M(t,x)$ and $h(t,x,w)=f(t,x,w-M(t,x))$. Then
\begin{equation}\label{trans}
 \partial_t h+\partial_x\big((w-M)h\big)=0 .
\end{equation}
At fixed $w$ the flow $\dot X=w-M(t,X)$ transports the initial distribution according to \eqref{trans}. The particle equation with these coefficients is renormalized. On the reflected lift $V=w-M(t,X)$ and $z=V+P(t,X)$ satisfy
\begin{equation}\label{position}
 z=w-rX,\qquad \dot z+rz=rP(t,X),\qquad
 |V(t)|\le e^{-rt}|v_0|+2\sup_{s\le t}\|P(s)\|_\infty .
\end{equation}
\end{lemma}
\begin{proof}
Take $\phi\in C_c^\infty([0,T)\times\R_x\times\R_w)$ on the covering line. Extend $M$ by $M(-t,x)=M(t,x)$ at zero and convolve in $(t,x)$, writing $M_\varepsilon$ for the result. On a compact strip containing the $(t,x)$-support of $\phi$, the extension is continuous and belongs to $W^{1,1}$ across zero. Consequently
\[
 M_\varepsilon\longrightarrow M\ \text{uniformly},\qquad
 (M_\varepsilon)_t\longrightarrow-m,\quad (M_\varepsilon)_x\longrightarrow\rho
       \quad\text{in }L^1\text{ on the original strip}.
\]
The smooth test $\phi(t,x,v+M_\varepsilon)$ is supported in a fixed velocity interval $|v|\le R$. Its derivatives give $\phi_t+v\phi_x$ and the additional term
\[
 \phi_w\big((M_\varepsilon)_t+m+v((M_\varepsilon)_x-\rho)\big),
\]
with the derivatives of $\phi$ evaluated at $(t,x,v+M_\varepsilon)$. After multiplication by $f$, the last term has integral bounded in absolute value by
\[
 2R\|f\|_\infty\|\phi_w\|_\infty
 \big(\|(M_\varepsilon)_t+m\|_1+R\|(M_\varepsilon)_x-\rho\|_1\big)\longrightarrow0.
\]
The terms $\phi_t+v\phi_x$ and the initial test converge uniformly on their common support. Changing variables $w=v+M(t,x)$ at fixed $(t,x)$ now gives \eqref{trans}, with initial density $h(0,x,w)=f_0(x,w-M(0,x))$.

Testing first by products in $(t,x)$ and $w$ gives the continuity equation at almost every fixed $w$. Its velocity $w-M(t,x)$ has at most linear growth and spatial Lipschitz constant $\|\rho(t)\|_\infty\in L^1(0,T)$. The characteristic formula transports the initial density and respects $(x,w)\mapsto(x+\ell,w+r\ell)$. For an admissible renormalization $\beta$, the spatial derivative of this velocity is $-\rho$, so the source in the equation for $\beta(h)$ is $\rho(h\beta'(h)-\beta(h))$. Testing back with $\psi(t,x,w-M_\varepsilon)$ uses the same $L^1$ derivative convergence, since $w-M$ is bounded on the test support. This gives \eqref{kinren}, including its initial value. Finally $z=w-rX$ gives $\dot z=-rV$ and hence \eqref{position}.
\end{proof}

For a drag $d(\rho)(u-v)$, the same differentiation along a smooth trajectory gives
\begin{equation}\label{dragtest}
 \frac{\dd}{\dd t}\big(V+\kappa M(t,X)\big)
   =\big(d(\rho)-\kappa\rho\big)(u-V).
\end{equation}
Within this class of additive labels, cancellation for all states and velocities is equivalent to $d(q)=\kappa q$. In particular, the conserved quantity in \eqref{twolabels} uses the linear dependence of the drag on $\rho$.

\begin{lemma}[Dependence of position on the particle label]\label{jaclemma}
Let $r>0$ and let $P,\rho$ be spatially periodic with $P_x=\rho-r$, $\rho\ge0$, $\|\rho(t)\|_\infty\in L^1(0,T)$ for all $T$, and $\sup_{t\ge1}\|\rho(t)\|_\infty\le B_*$. Assume $P$ is continuous and
$\|P(t)\|_\infty+\|\rho(t)-r\|_1\le Ce^{-\beta t}$ for $t\ge1$.
At fixed $a\in[0,\ell)$ solve $\dot X_t(y)=r(y-X_t(y))-P(t,X_t(y))$, $X_0(y)=a$. For $t\ge2$ the map $X_t$ is increasing and bi-Lipschitz, and its inverse has a Lipschitz constant independent of $a,t$. For each bounded interval $K$ and $0<\sigma<\min(r,\beta)$,
\begin{equation}\label{jacclose}
 \sup_{y\in K}|X_t(y)-y|+\int_K|X'_t(y)-1|\dd y\le C_K e^{-\sigma t}.
\end{equation}
The prime denotes differentiation in $y$, and $C_K$ is uniform in $a$.
\end{lemma}
\begin{proof}
Subtract the trajectory equations at $y_2>y_1$. The difference $d(t)=X_t(y_2)-X_t(y_1)$ satisfies
\[
 d'+a(t)d=r(y_2-y_1),\qquad d(0)=0,
\]
where $a(t)$ is a difference quotient of $M=rx+P$ between the two positions (zero when they coincide). Since $M_x=\rho$, we have $0\le a(t)\le\|\rho(t)\|_\infty$. Integrating, and retaining the last unit of time for the lower bound, yields
\begin{equation}\label{bilip}
 c_0(y_2-y_1)\le X_t(y_2)-X_t(y_1)\le rt(y_2-y_1)\qquad(t\ge2).
\end{equation}
Here $c_0=r\int_0^1 e^{-B_*s}\dd s>0$. Integrating from zero also gives $d(t)\ge rt\,e^{-\int_0^t\|\rho(s)\|_\infty\dd s}(y_2-y_1)$ for every $t>0$. Hence the inverse map takes null sets to null sets. Differentiating the integral equation by the Lipschitz chain rule gives, for a.e. $y,t$,
\begin{equation}\label{jac}
 \partial_t X'_t=r-\rho(t,X_t)X'_t,\qquad X'_0=0 .
\end{equation}
The bound $0\le X'_t\le rt$ justifies differentiating the integral in time by Fubini. Integrating $\partial_t(X_t-y)+r(X_t-y)=-P(t,X_t)$ gives the uniform position estimate in \eqref{jacclose}. To estimate the Jacobian, subtract its equilibrium value one in \eqref{jac} and use $x=X_t(y)$. This gives
\begin{align*}
 (\partial_t+r)(X'_t-1)&=-(\rho(t,X_t)-r)X'_t,\\
 \int_K|\rho(t,X_t)-r|X'_t\dd y
   &=\int_{X_t(K)}|\rho(t,x)-r|\dd x\le C_K\|\rho(t)-r\|_1 .
\end{align*}
The supremum bound already proved places $X_t(K)$ in a fixed bounded interval. Periodicity gives the last inequality. With $A_K(t)=\int_K|X'_t-1|$, integration yields
\[
 A_K(t)\le e^{-r(t-2)}A_K(2)+C_K\int_2^t e^{-r(t-s)}\|\rho(s)-r\|_1\dd s ,
\]
which proves \eqref{jacclose}, including the case $r=\beta$ with the stated strict choice of $\sigma$.
\end{proof}

\begin{lemma}[Testing in a transported coordinate]\label{weightlemma}
On $[0,L]$ with zero boundary flux, or on $\T_\ell$, let a nonnegative density $\mu$ be bounded on each finite time interval and satisfy
\[
 \partial_t\mu+\partial_x(u\mu)=g(t)\mu,\qquad \mu(0)=\mu_0
\]
in the weak sense, with $g\in L^1_{\rm loc}([0,\infty))$ and $u\in L^2_{\rm loc}([0,\infty);L^\infty)$. Suppose $S$ is continuous, $S_t+uS_x=0$, $S_t\in L^1_{\rm loc}$ and $S_x\in L^2_{\rm loc}$ up to $t=0$. Then for every continuous test function $\varphi$ on the range of $S$,
\begin{equation}\label{transportidentity}
 \int\varphi(S(t,x))\mu(t,x)\dd x
  =\exp\!\left(\int_0^t g(s)\dd s\right)
       \int\varphi(S(0,a))\mu_0(a)\dd a
\end{equation}
for a.e. $t$. On the torus use periodic tests and a lift of the coordinate.
\end{lemma}
\begin{proof}
Test by $\eta(t)\varphi(S(t,x))$ with $\eta$ smooth and compactly supported in time. For smooth $\varphi$, mollify $S$ up to $t=0$. The convergence is uniform, with time derivatives converging in $L^1$ and spatial derivatives in $L^2$, as required by $\mu\in L^\infty$ and $u\in L^2_tL^\infty_x$ on each finite strip. Since $S_t+uS_x=0$, the integral in \eqref{transportidentity} solves $Y'=gY$ with the indicated initial value. Uniform approximation by a countable dense family of smooth tests gives the assertion.

\end{proof}

\section{Convergence of the particle density}\label{particles}
We fix a solution from Appendix~\ref{construction} with $r>0$. Unless the interval is specified, integrals in this section are over the reflected torus and $v\in\R$. With $M_f=\iint f_0$ put
\begin{align}
 E(t)&=\int\left(\frac{|m|^2}{2\rho}+\frac{\rho^\gamma}{\gamma-1}\right)\dd x
           +\tfrac12\iint v^2 f\dd x\dd v,\label{energy}\\
  D(t)&=\nu\int|u_x|^2\dd x+\iint\rho f|u-v|^2\dd x\dd v,\qquad
 \E(t)=E(t)-\frac{\ell r^\gamma}{\gamma-1}.\label{diss}
\end{align}
On $I$ the same definitions use interval integrals and $L$ in place of $\ell$.

\subsection{Decay of the relative energy}
For $t\ge1$ the density is bounded above by a constant $B_*$ by \eqref{densitybound}. On the interval $0\le q\le B_*$ the relative pressure satisfies
\begin{equation}\label{comp}
 c(q-r)^2\le q^\gamma-r^\gamma-\gamma r^{\gamma-1}(q-r)
             \le C(q-r)^2 .
\end{equation}
Indeed the quotient by $(q-r)^2$ extends to a strictly positive continuous function at $q=r$ and is positive at $q=0$. Also $(q^\gamma-r^\gamma)(q-r)\ge c(q-r)^2$. We will show that for $t\ge2$,
\begin{equation}\label{exp}
 \E(t)+\int_t^\infty D(s)\dd s\le Ce^{-ct},\qquad
 \|\rho(t)-r\|_2+\|P(t)\|_\infty\le Ce^{-ct/2}.
\end{equation}

The energy equality gives $\E'=-D$. To obtain a dissipative term for $\rho-r$, test the momentum equation with the periodic primitive $P$ of $\rho-r$, as in Fang, Zi and Zhang \cite{FZZ} for the compressible fluid and Li and Shou \cite[(5.18)]{LS} for this coupling. Set $J(t)=\int mP$, $p=\rho^\gamma$ and $F=\rho\int_\R (u-v)f\dd v$. Fix $\eta\in C_c^\infty(1,\infty)$ and test momentum with $\eta P_\varepsilon$, where $P_\varepsilon$ is a space-time convolution, periodic in $x$. On a slightly larger compact time interval, $P$ is continuous, $P_t=-m\in L^2_{t,x}$ and $P_x=\rho-r$ is bounded by \eqref{densitybound}. Thus $P_\varepsilon\to P$ uniformly, $(P_\varepsilon)_t\to-m$ in $L^2$, and $(P_\varepsilon)_x\to\rho-r$ almost everywhere with a common bound. The time derivative pairs with $m\in L^2$. The spatial derivative pairs with $mu+p-\nu u_x\in L^1$. The uniform convergence suffices for the drag $F\in L^1$. Passing to the limit gives
\[
 -\int\eta' J\dd t=\int\eta\!\int\big[mP_t+(mu+p-\nu u_x)P_x-FP\big]\dd x\dd t .
\]
The right side is locally integrable in time. The weak continuity of $m$ and uniform continuity of $P$ identify $J$ with its absolutely continuous representative on compact subintervals. Now $mP_t+muP_x=-r\rho u^2$ and $\int(\rho-r)=0$, which give
\begin{align}
 J'&=\int (p-r^\gamma)(\rho-r)-\nu\int u_x(\rho-r)
                -r\int\rho u^2-\int FP \notag\\
   &\ge c_1\|\rho-r\|_2^2-C_1D. \label{J}
\end{align}
Here $\|P\|_\infty\le C_L\|\rho-r\|_2$ and oddness yields $\|u\|_\infty\le C_L\|u_x\|_2$. For the force term use Cauchy--Schwarz in the form $\|F\|_1^2\le B_*M_f\iint\rho f|u-v|^2$, followed by Young's inequality.

The drag dissipation is weighted by $\rho$, so it gives no direct control of particle velocities in a vacuum region. The conserved quantity in Lemma~\ref{masslemma} suggests using $z=v+P$ instead. Its second moment $K(t)=\iint(v+P)^2f$ satisfies
\begin{equation}\label{Zdiff}
 K'=-2rK+2r\iint(v+P)Pf\le-rK+C_2\|\rho-r\|_2^2 .
\end{equation}
This follows by integrating $\frac{\dd}{\dd t}z^2=-2rz^2+2rzP(t,X)$ against the initial particle measure. The second moment bound and \eqref{position} justify the differentiation.

Take $\mathcal L=\E-\delta J+\delta^2 K$. We have $|J|+K\le C\E$ by \eqref{comp} and conservation of mass. Choose $\delta>0$ small. In \eqref{J} the error $\delta C_1D$ is absorbed by the energy dissipation, and in \eqref{Zdiff} the error $\delta^2 C_2\|\rho-r\|_2^2$ is absorbed by $\delta c_1\|\rho-r\|_2^2$. Thus
\[
 \tfrac12\E(t)\le\mathcal L(t)\le2\E(t),\qquad
 \mathcal L'+c\big(D+\|\rho-r\|_2^2+K\big)\le0 .
\]

To bound the relative energy by the quantities in parentheses, use \eqref{comp} for the pressure, Poincar\'e's inequality for the fluid kinetic energy, and $\iint v^2f\le2K+2M_f\|P\|_\infty^2$ for the particles. Gronwall's inequality gives \eqref{exp}, since \eqref{allenergy} then gives $\int_t^\infty D(s)\dd s=\E(t)$.

\subsection{The limiting particle position}

For a particle issued from $(a,v_0)$ write $y=S_0(a)+v_0/r$. Lemma~\ref{masslemma} gives $r(y-X)=z$. Together with \eqref{exp}, the damped equation \eqref{position} gives $X(t)\to y$ and $V(t)\to0$. For $M_f>0$, couple each trajectory with its terminal point $([y]_\ell,0)$ using the same initial particle. Then
\begin{equation}\label{w2}
 W_2^2\!\left(\frac{f(t)\dd x\dd v}{M_f},
       \frac{n_\infty(x)\dd x\otimes\delta_0}{M_f}\right)
 \le \frac{r^{-2} K(t)+\iint v^2f(t)}{M_f}\le Ce^{-ct}.
\end{equation}
Here $W_2$ is the quadratic Wasserstein distance, with cost equal to squared torus distance plus squared velocity difference. Changing variables $v_0=r(y-S_0(a))$ at fixed $a$ gives on the torus
\begin{equation}\label{profile}
 n_\infty(x)=r\int_0^\ell\sum_{j\in\mathbb Z} f_0\big(a,r(x+j\ell-S_0(a))\big)\dd a .
\end{equation}
Folding positions contracts distances and reversing velocity preserves its distance to zero. Both reflected copies have the same folded measure. Dividing by two yields \eqref{physprofile} and \eqref{w2} with the interval mass. If $M_f=0$ both particle densities vanish.

\subsection{Strong convergence and sharpness}\label{marginal}
\begin{proposition}\label{relax}
Assume $s\ge2$ and $\iint |v|^s f_0<\infty$. Then
\begin{equation}\label{strongL3}
 \sup_{t\ge2}\|n(t)\|_{s+1}+\|n_\infty\|_{s+1}<\infty,\qquad
 n(t)\to n_\infty\quad\text{in }L^{s+1}.
\end{equation}
For bounded initial velocity support, $\sup_{t\ge2}\|n(t)\|_\infty<\infty$, convergence holds in every finite $L^p$, and
\begin{equation}\label{support}
 \operatorname{ess\,supp}_v f(t)\subset[-Ce^{-\sigma t},Ce^{-\sigma t}],\qquad t\ge2 .
\end{equation}
With bounded initial velocity support and $\int\TV_v(f_0(a,\cdot))\dd a<\infty$, the convergence is exponential in every finite $L^p$. Here $\TV_v$ denotes total variation in velocity.
\end{proposition}
\begin{proof}
Fix the initial position $a$. In the coordinate $y=S_0(a)+v_0/r$ the initial density is $g(y)=r f_0(a,r(y-S_0(a)))$. At time $t$ it is carried to $X_t(y)$ and then reduced modulo $\ell$. We first compare the density at time $t$ with $g$ for bounded labels. For compactly supported $g$ of bounded variation, choose $K$ containing both its support and its inverse image under $X_t$. For a common bounded velocity support, \eqref{bilip} and \eqref{jacclose} allow $K$ to be fixed uniformly in $a$ and $t\ge2$. With $h=\sup_K|X_t(y)-y|$, change of variables gives the following bound on the real-line $L^1$ difference.
\[
 \int_{\R}|g(y)-g(X_t(y))X'_t(y)|\dd y
 \le 2h\,\TV(g)+\|g\|_\infty\int_K |X'_t(y)-1|\dd y .
\]
The first term follows by integrating the variation of $g$ over $[y-h,y+h]$ and applying Fubini. Both terms decay exponentially by Lemma~\ref{jaclemma}. Integrating in $a$ and projecting to the torus preserves the bound. Approximation in $L^1(\dd a\dd y)$ by bounded functions compactly supported and smooth in $y$ gives $L^1$ convergence for general data, since transport and projection contract total variation. For bounded initial velocities, $V=r(y-X)-P$ also gives \eqref{support}.

To reach the endpoint $s+1$ one must control large initial velocities uniformly in time. The contribution from $|v_0|\le H$, $H\ge1$, is bounded pointwise by $CH\|f_0\|_\infty$. Indeed at fixed $a$ these labels and their positions lie in intervals of length at most $CH$ by \eqref{position}. At most $CH$ preimages project to a single point of the torus, and \eqref{bilip} bounds the inverse Jacobian. For $R>1$ let $n_R(t,x)$ denote the density from $|v_0|>R$, and let $e_R(t,x)$ denote the same density with initial weight $|v_0|^s$. Splitting at $H$ gives
\[
 n_R\le CH\|f_0\|_\infty+H^{-s}e_R,\qquad n_R\le e_R .
\]
Where $e_R\ge1$ choose $H=e_R^{1/(s+1)}$ to balance the two terms. Where $e_R<1$ use $n_R\le e_R$. Integration in $x$ yields
\begin{equation}\label{endpoint}
 \|n_R(t)\|_{s+1}^{s+1}\le C_s\int e_R(t,x)\dd x
   =C_s\iint_{|v_0|>R}|v_0|^s f_0(a,v_0)\dd a\dd v_0 .
\end{equation}
For the limiting position $X_\infty(y)=y$ the same count gives the same bound. At fixed $R$, the contribution from $|v_0|\le R$ converges in $L^{s+1}$ by its uniform $L^\infty$ bound and the $L^1$ convergence. The right side of \eqref{endpoint} tends to zero with $R\to\infty$. This proves \eqref{strongL3} and Theorem~\ref{main}\,\textup{(i)}. For bounded velocity support and integrable velocity variation, interpolation with the exponential $L^1$ estimate gives the asserted rates in every finite $L^p$.
\end{proof}

\begin{proposition}[Sharpness]\label{optimalendpoint}
For every $s\ge2$ there are data satisfying \eqref{data} with a finite $s$-th moment for which the limiting marginal belongs to $L^{s+1}$ and to no $L^p$ with $p>s+1$. The example can be chosen with constant initial fluid density and specular symmetry.
\end{proposition}
\begin{proof}
Formula \eqref{profile} adds the contributions of velocities whose labels differ by one period. With $\ell=r=1$, placing unit density at $H(y)$ such velocities gives a limiting density of order $H(y)$ at a moment cost of order $H(y)^{s+1}$. Take $\rho_0=1$, $m_0=0$, and $a\in[0,1)$. For $0<y<1/4$ choose
\[
 H(y)=\left\lfloor y^{-1/(s+1)}(\log(e/y))^{-2/(s+1)}\right\rfloor ,
\]
and set $H=0$ elsewhere on $[0,1)$. Each pair $(a,v)$ gives unique $k\in\mathbb Z$ and $y\in[0,1)$ with $a+v=k+y$. Define $f^+(a,v)=\mathbf1_{\{1\le k\le H(y)\}}$ on the fundamental interval and extend periodically in $a$. Take $f_0(a,v)=f^+(a,v)+f^+([-a]_1,-v)$. The distribution is bounded and symmetric. Its mass and $s$-th moment are bounded by $C_s\int_0^1(H+H^{s+1})\dd y<\infty$, by summation over $k$. By \eqref{profile} its limiting marginal is $H(y)+H([-y]_1)$ almost everywhere. The integral of $H^p$ diverges at zero for every $p>s+1$, whereas $H^{s+1}$ is integrable by the logarithmic factor. The density is even, so restriction to the physical half-period retains the divergence.
\end{proof}

Boundedness of $f_0$ matters even with bounded velocities. For $I=(0,1)$ and $\rho_0=1$, let $f_0(a,v)=g(a+v)$ with $g(y)=y^{-1/2}$ on $(0,1)$ and $g=0$ elsewhere. Then $|v|<1$ on its support, and \eqref{physprofile} gives the density $g(x)$, which does not belong to $L^3(I)$.

\begin{corollary}[Stability of the terminal density]\label{stable}
Fix $s\ge2$. On the same interval $I$, let the limiting data, of mean density $r>0$, and a sequence of data satisfy the hypotheses of Proposition~\ref{relax}. Assume
\[
 \rho_{0,k}\to\rho_0\ \text{in }L^1(I),\qquad
 f_{0,k}\to f_0\ \text{in }L^1(1+|v|^s),\qquad
 \sup_k\|f_{0,k}\|_\infty<\infty .
\]
Then $n_{\infty,k}\to n_\infty$ in $L^{s+1}(I)$.
\end{corollary}
\begin{proof}
We have $r_k\to r>0$ and $S_{0,k}\to S_0$ uniformly. First restrict to $|v|\le H$. Replacing $f_{0,k}$ by $f_0$ on this set costs at most $\|f_{0,k}-f_0\|_1$ in total variation. For fixed $a$ the label is $S_{0,k}(a)+v/r_k$. Continuity of translations and dilations in $L^1_v$, followed by dominated convergence in $a$ and reflection, gives $L^1$ convergence of the truncated marginals.

The preimages of $0<x<L$ under reflection are $2jL\pm x$. At fixed $a$ each sign allows at most $1+H/(Lr_k)$ such positions. The truncated densities are therefore bounded by $2(Lr_k+H)\|f_{0,k}\|_\infty$, giving convergence in $L^{s+1}$ by interpolation. The same count in \eqref{endpoint} bounds the terminal tails with a constant depending only on $L,s,\sup r_k,\sup\|f_{0,k}\|_\infty$. Weighted $L^1$ convergence makes these tails uniformly small as $H\to\infty$.

\end{proof}

\section{The limiting vacuum measure}\label{vacuum}
Set $k=r^\gamma/\nu$ with $r>0$. We work on $[0,L]$ since vacuum may concentrate at either wall. For the uncoupled fluid, questions of separation of vacuum and of finite-time closure for density-dependent viscosity are treated in \cite{DuanLuo,LGL}.

\begin{proposition}\label{vac}
There are constants $c,C,\sigma>0$ and a measure $\lambda_\infty$ on $[0,L]$ such that
\begin{equation}\label{vacrate}
 c|Z_0|e^{-kt}\le |Z_t|\le C|Z_0|e^{-kt}\qquad\text{for a.e. }t>0.
\end{equation}
For every Lipschitz function $\varphi$ on $[0,L]$ and a.e. $t\ge2$,
\begin{equation}\label{vacbl}
 \left|e^{kt}\int_{Z_t}\varphi(x)\dd x-\int_{[0,L]}\varphi\dd\lambda_\infty\right|
 \le C e^{-\sigma t}\big(\|\varphi\|_\infty+\operatorname{Lip}\varphi\big).
\end{equation}
For every Borel set $A\subset[0,L]$, its limiting mass satisfies
\begin{equation}\label{vaclabel}
 c\,|Z_0\cap S_0^{-1}(A)|\le\lambda_\infty(A)\le C\,|Z_0\cap S_0^{-1}(A)|.
\end{equation}
\end{proposition}

\begin{proof}
Write $\chi=\mathbf1_{Z_t}$. In the renormalized continuity equation use $\beta_\varepsilon(q)-1$, where $\beta_\varepsilon(q)=\varepsilon/(\varepsilon+q)$, and add the constant balance. Both $\beta_\varepsilon$ and $\beta_\varepsilon-q\beta_\varepsilon'$ are uniformly bounded and tend to $\mathbf1_{\{q=0\}}$. Dominated convergence up to $t=0$ gives
\begin{equation}\label{chi}
 \chi_t+(u\chi)_x=u_x\chi,\qquad \chi(0)=\mathbf1_{Z_0}.
\end{equation}
The particle mass to the left of $x$ supplies the correction to the fluid momentum. Work on the physical interval and set
\begin{equation}\label{vacprimitive}
 N(t,x)=\int_0^x n(t,y)\dd y,\qquad q=m+\rho N,\qquad
 b_x=q,\quad \int_0^L b(t,x)\dd x=0 .
\end{equation}
At the smooth approximation level, specular reflection gives $j(t,0)=0$ for $j=\int_\R vf\dd v$. Integrating $n_t+j_x=0$ from the wall gives $N_t=-j$. Hence $N_t+uN_x=un-j$, and the balance for $\rho N$ cancels the drag in momentum. It follows that
\begin{align}
 q_t+(uq+p)_x&=\nu u_{xx},\qquad p=\rho^\gamma,\notag\\
 b_t+ub_x&=\nu u_x-p+A(t),\qquad A(t)=\frac1L\int_0^L(uq+p)\dd x .\label{scalarpressure}
\end{align}
Indeed the spatial derivative of $b_t+uq+p-\nu u_x$ vanishes. Taking its mean gives $A(t)$ since $b$ has mean zero and $u$ vanishes at both walls. This avoids any boundary trace of $u_x$ or $p$ for weak solutions. Appendix~\ref{weakprimitive} passes this identity from the approximations, with $b$ continuous and bounded, $b_t\in L^1$ and $b_x\in L^2$ up to $t=0$.

On $Z_t$ the pressure vanishes. Multiplying \eqref{chi} by $e^{kt-b/\nu}$ now cancels $u_x\chi$ and gives
\begin{equation}\label{K}
 \mu=e^{kt-b/\nu}\chi,\qquad
 \mu_t+(u\mu)_x=\frac{r^\gamma-A(t)}{\nu}\,\mu .
\end{equation}
Mollification of $b$ justifies the product up to time zero by the preceding derivative bounds and $u\in L^2_tL^\infty_x$. Reflection gives the zero-flux identity on $[0,L]$.

Since $0\le N\le\int_I n_0$, we have
\[
 |A(t)-r^\gamma|\le C\big(\E(t)+\sqrt{\E(t)}\big).
\]
Indeed $\int_I|uq|\le\int_I\rho u^2+\|N\|_\infty\int_I\rho|u|$, and $L^{-1}\int_I p-r^\gamma$ is controlled by the relative internal energy. Thus
\[
 \omega(t)=\exp\!\left(\nu^{-1}\int_0^t(r^\gamma-A(s))\dd s\right)
\]
has a positive limit $\omega_\infty$, with $|\omega(t)-\omega_\infty|\le Ce^{-\sigma t}$ by \eqref{exp}. Lemma~\ref{weightlemma} applies to \eqref{K} with $S=M/r$, since $S_t=-m/r$, $S_x=\rho/r$ belong to $L^2$ on finite strips by \eqref{pressuregain}. It gives the finite-time identity
\begin{equation}\label{vacidentity}
 \int_0^L \varphi(S(t,x))\mu(t,x)\dd x
 =\omega(t)\int_{Z_0}\varphi(S_0(a)) e^{-b(0,a)/\nu}\dd a .
\end{equation}
Taking $\varphi=1$ and using boundedness of $b$ proves \eqref{vacrate}.

We already know that $S(t,x)-x=P(t,x)/r$ tends to zero exponentially. To find the limit of $b$, put $H(t,x)=r\int_0^L(x-y)_+ n(t,y)\dd y$, where $s_+=\max(s,0)$. Since $H_x=rN$,
\[
 \|b-(H-L^{-1}\textstyle\int_I H)\|_\infty
   \le C\big(\|m\|_1+\|N\|_\infty\|\rho-r\|_1\big).
\]
The right side decays exponentially by \eqref{exp}. For every $x$, the function $y\mapsto(x-y)_+$ has Lipschitz constant at most one. The interval coupling in \eqref{w2} therefore gives uniform exponential convergence of $H(t,\cdot)$ to $H_\infty(x)=r\int_0^L(x-y)_+n_\infty(y)\dd y$ (both vanish in the absence of particles). Thus $b(t,\cdot)$ converges at the same rate to the Lipschitz function $b_\infty=H_\infty-L^{-1}\int_I H_\infty$.

Use \eqref{vacidentity} with the test $e^{b_\infty/\nu}\varphi$. The mass $\int\mu(t)$ is bounded, and both $\|S(t)-x\|_\infty$ and $\|b(t)-b_\infty\|_\infty$ are at most $Ce^{-\sigma t}$. Thus replacing $e^{b_\infty(S(t,x))/\nu}\varphi(S(t,x))$ by $e^{b(t,x)/\nu}\varphi(x)$ costs at most the right side of \eqref{vacbl}. This proves \eqref{vacbl} with
\begin{equation}\label{vaclimit-test}
 \int_{[0,L]}\varphi\dd\lambda_\infty
  =\omega_\infty\int_{Z_0}\varphi(S_0(a))W_0(a)\dd a,\qquad
 W_0(a)=\exp\!\left(\frac{b_\infty(S_0(a))-b(0,a)}{\nu}\right).
\end{equation}
This identity determines a finite Borel measure. Its weight is bounded above and below by positive constants, proving \eqref{vaclabel}. Since $S_0$ is absolutely continuous and $S_0'=\rho_0/r$ vanishes on $Z_0$, its image $S_0(Z_0)$ has length zero by the one-dimensional area formula. This proves singularity. Testing by one gives \eqref{vaclim}, and uniform approximation gives Theorem~\ref{main}\,\textup{(ii)} for every continuous test.
\end{proof}

\begin{corollary}[Relative weights of the vacuum set]\label{vacshape}
Let $r>0$ and $|Z_0|>0$. For any Borel set $A\subset[0,L]$ the fraction of the limiting vacuum mass in $A$ is
\begin{equation}\label{vacnormalized}
 \frac{\lambda_\infty(A)}{V_*}
  =\frac{\displaystyle\int_{Z_0\cap S_0^{-1}(A)} W_0(a)\dd a}
         {\displaystyle\int_{Z_0} W_0(a)\dd a}.
\end{equation}
Here $S_0^{-1}(A)=\{a:S_0(a)\in A\}$. Up to an additive constant independent of $a$, the logarithm of the weight is given by
\begin{equation}\label{elementaryweight}
 \nu\log W_0(a)=r\int_0^{S_0(a)}(S_0(a)-y)n_\infty(y)\dd y
       -\int_0^a\big(m_0(x)+\rho_0(x)N_0(x)\big)\dd x
       +\text{\rm constant},
\end{equation}
where $N_0(x)=\int_0^x n_0(y)\dd y$ and $n_0=\int_\R f_0\dd v$.
\end{corollary}
\begin{proof}
Apply the measure identity \eqref{vaclimit-test} to $A$ and to $[0,L]$ and divide. Formula \eqref{elementaryweight} follows from $b_x(0)=m_0+\rho_0N_0$ and the formula for $b_\infty$. The additive constant comes from the zero means of the two primitives and cancels in \eqref{vacnormalized}. The terminal density $n_\infty$ is given by \eqref{physprofile}. Thus data admissible at two pressure exponents give the same normalized vacuum measure at fixed viscosity. The absolute mass $V_*=\omega_\infty\int_{Z_0}W_0(a)\dd a$ still involves $\omega_\infty$ from the evolution, with the mean-zero normalization in \eqref{vaclimit-test}.
\end{proof}

\begin{remark}
Consider two gaps of equal length. On $(0,1)$ let $\rho_0$ take the values $1,0,1,0,3$ on the five successive intervals of length $1/5$. Take $m_0=0$ and $f_0(a,v)=(N/2)\mathbf1_{\{|v|<1\}}$ with $N>0$. Since $r=1$ and the velocities fill one unfolded period uniformly, $n_0=n_\infty=N$. Label the two gaps from left to right. Their mass coordinates are $1/5$ and $2/5$, and the weight is constant on each gap. The first integral on the right of \eqref{elementaryweight} increases by $\frac{N}{2}((2/5)^2-(1/5)^2)=3N/50$. The second increases by $\int_{2/5}^{3/5}Nx\dd x=N/10$. Their difference is $-N/25$. The gaps have equal length, so \eqref{vacnormalized} gives
\[
 \frac{\lambda_\infty(\{2/5\})}{\lambda_\infty(\{1/5\})}=e^{-N/(25\nu)} .
\]
A singular continuous limit is also possible. Let $\rho_0$ vanish on a closed nowhere dense set of positive length and equal a positive constant elsewhere. Take $m_0=0$ and bounded particle data with compact velocity support. Every open interval contains fluid mass, so $S_0$ is strictly increasing. Each point has at most one preimage, and \eqref{vacnormalized} gives zero mass to every singleton.

If $|Z_0|>0$, the vacuum contributes at least $r|Z_t|^{1/p}$ to $\|\rho(t)-r\|_{L^p(I)}$ for $1\le p<\infty$ and gives $\|\rho(t)-r\|_\infty\ge r$ at almost every time. The relative internal energy density equals $r^\gamma$ on $Z_t$. Thus $\E(t)\ge c e^{-r^\gamma t/\nu}$. The optimal decay rate of the full relative energy remains open.
\end{remark}

\appendix
\section{Construction and energy equality}\label{construction}
\subsection{Weak formulation and approximation}
We prove the finite-time results used in the main text for every fluid mass, including zero. The primitive $P=\I[\rho]$ exists for all such masses and equals $r(S-x)$ in the positive-mass case. On the reflected lift $\I$ denotes the zero-mean periodic primitive, $M_f=\iint f_0$ the particle mass and $j=\int_\R vf\dd v$ the particle momentum density.

We use the vacuum convention in \eqref{data} for $|m|^2/\rho$. We evaluate $E$ from \eqref{energy} using the time representatives of $\rho,m,f$. The dissipation \eqref{diss} is defined almost everywhere. A finite-energy weak solution has, on each finite time interval $(0,T)$,
\[
 \rho\in L^\infty_t L^\gamma_x,\quad \sqrt{\rho}u\in L^\infty_t L^2_x,
 \quad u\in L^2_t H^1_0(I),\quad f\ge0,
\]
\[
 (1+v^2)f\in L^\infty_t L^1_{x,v},\quad
 \rho f(u-v)\in L^1_{t,x,v},\quad  D\in L^1(0,T).
\] The first two equations of \eqref{sys} hold distributionally with initial data. Continuity admits tests on $\overline I$ and momentum tests vanish at both endpoints. The kinetic equation reads
\begin{equation}\label{spec}
 \int_0^\infty\iint_{I\times\R}
 f\big(\psi_t+v\psi_x+\rho(u-v)\psi_v\big)\dd x\dd v\dd t
 +\iint_{I\times\R} f_0\psi(0)\dd x\dd v=0
\end{equation}
for compactly time- and velocity-supported $C^1$ tests on the closed interval satisfying $\psi(t,x,v)=\psi(t,x,-v)$ at $x=0,L$. 

\begin{theorem}[Global existence and kinetic compactness]\label{exist}
Under \eqref{data}, there exists a global finite-energy weak solution of \eqref{sys}, satisfying \eqref{spec}, the renormalized continuity equation and
\begin{equation}\label{densitybound}
 \|\rho(t)\|_\infty\le C(1+t^{-1/\gamma}),\quad t>0,\qquad
 f\in L^\infty((0,T)\times I\times\R)\quad\text{for every }T<\infty.
\end{equation}
The constant in the density bound depends only on the initial energy, the two masses, $\nu,\gamma,L$. For bounded $\rho_0$ a bound uniform down to $t=0$ may also depend on $\|\rho_0\|_\infty$. The density, momentum and particle distribution are weakly continuous in time with values in $L^\gamma$, $L^{2\gamma/(\gamma+1)}$ and $L^1$, respectively. Both masses are conserved, and the total energy of these representatives satisfies
\begin{equation}\label{allenergy}
 E(t)+\int_s^t D(\tau)\dd\tau=E(s),\qquad 0\le s\le t .
\end{equation}
For each $\beta\in C^1([0,\infty))$ with $\beta(0)=0$ and $\beta'$ bounded on the range of $f$ on each finite time interval, the renormalized equation reads
\begin{equation}\label{kinren}
 \partial_t\beta(f)+v\partial_x\beta(f)
 +\partial_v\big(\rho(u-v)\beta(f)\big)=\rho\big(f\beta'(f)-\beta(f)\big).
\end{equation}
A sequence of positive-density approximations can be chosen with
\begin{equation}\label{strongf}
 f_k\longrightarrow f\quad\text{strongly in }
 L^1((0,T)\times I\times\R;(1+v^2)\dd t\dd x\dd v).
\end{equation}
\end{theorem}
The particle evolution is unique for a prescribed fluid. Uniqueness for the coupled system in the full finite-energy class remains open.

\smallskip\noindent\emph{Positive-density approximation.}
Extend $\rho_0$ evenly, $m_0$ oddly and $f_0$ by $f_0(-x,-v)=f_0(x,v)$ to $\T_\ell$. In this construction we allow $r=0$ and always take $r_k=\av{\rho_{0,k}}>0$.

Choose smooth symmetric data $(\rho_{0,k},u_{0,k},f_{0,k})$, with $\rho_{0,k}>0$ and $f_{0,k}\ge0$ compactly supported in velocity, such that
\begin{equation}\label{approx}
 \rho_{0,k}\to\rho_0\text{ in }L^\gamma,\qquad
 \sqrt{\rho_{0,k}}u_{0,k}\to m_0/\sqrt{\rho_0}\text{ in }L^2,\qquad
 f_{0,k}\to f_0\text{ in }L^1(1+v^2),
\end{equation}
with $\|f_{0,k}\|_\infty\le\|f_0\|_\infty$. For the fluid, approximate a truncated $u_0$ in $L^2(\rho_0\dd x)$ by smooth odd functions, then mollify the even density and add a positive constant, using a diagonal sequence. For the kinetic data use velocity truncation and symmetric nonnegative mollification. In particular $m_{0,k}\to m_0$ in $L^{2\gamma/(\gamma+1)}$ and $E_k(0)\to E(0)$.

Theorem 1.1 and (1.9)--(1.10) of \cite{LS} give strong solutions with $\rho_k,u_k\in C_tH^2_x$ and $f_k\in C_tC^1_{x,v}$. Indeed, rescaling $y=x/\ell$, $w=v/\ell$ changes $(\rho,u,f)$ into $(\rho(t,\ell y),u(t,\ell y)/\ell,\ell f(t,\ell y,\ell w))$. In \cite[(1.3)--(1.4)]{LS} this corresponds to $A=\ell^{-2}$, $\kappa=1$, $\beta=0$, and $\mu_0,\mu_1>0$ with $\mu_0+\mu_1=\nu/\ell^2$. Their formula (2.14) includes this constant-viscosity case. Related results concern classical solutions with reflecting boundaries \cite{Jiang}, a hydrodynamic limit \cite{CWY}, and weak solutions with velocity diffusion \cite{LSFP}. For the latter, convergence to a global Maxwellian is obtained under stronger assumptions. Uniqueness of each strong approximation preserves reflection symmetry, hence

\begin{equation}\label{oddmeans}
 \av{\rho_k u_k}=\av{j_k}=0,\qquad \av{u_k}=0.
\end{equation}
The basic energy identity and mass conservation are
\begin{equation}\label{base}
 E_k(t)+\int_0^t D_k(\tau)\dd\tau=E_k(0),\qquad
 \int\rho_k=\int\rho_{0,k},\quad\iint f_k=\iint f_{0,k}.
\end{equation}
By Poincar\'e and one-dimensional Sobolev embedding,
$\int_0^\infty(\|u_k\|_{H^1}^2+\|u_k\|_\infty^2)\dd t\le C$.

\subsection{Bounds on the fluid density and the drag}\label{initiallayer}
Restrict a positive-density approximation to $I$ and suppress its index. The correction $u+N$ in \eqref{vacprimitive} is the wall-normalized form of the effective velocity in \cite[(1.18)]{LS}. With $q=\rho(u+N)$, the smooth calculation in \eqref{scalarpressure} gives
\begin{equation}\label{logODE}
 (\partial_t+u\partial_x)(\nu\log\rho+b)+\rho^\gamma=A(t).
\end{equation}
We have $\|N\|_\infty\le\int_I n_0$ and $\|m\|_1^2\le(\int_I\rho)(\int_I\rho u^2)$. Thus $\|b\|_\infty\le C\|q\|_1$, and
\[
 L|A(t)|\le\int_I(\rho u^2+p)\dd x+\|N\|_\infty\|m\|_1.
\]
The energy and mass bounds in \eqref{base} control both quantities. Along the fluid flow, $y=\rho e^{b/\nu}$ therefore obeys $y'\le C y-cy^{\gamma+1}$ with $C,c>0$ independent of the approximation. Integration gives
\[
 y(t)^{-\gamma}\ge e^{-\gamma Ct}y(0)^{-\gamma}
              +\frac cC(1-e^{-\gamma Ct}).
\]
The second term alone yields
\begin{equation}\label{rhok}
 \|\rho_k(t)\|_\infty\le \rho_*(t),\qquad \rho_*(t)=C(1+t^{-1/\gamma}).
\end{equation}
By reflection this holds on the whole torus.

This bound is independent of the initial density bounds. Retaining the initial term gives a uniform estimate down to $t=0$ for bounded $\rho_0$. Since $\rho_*\in L^1(0,T)$ for $\gamma>1$, kinetic characteristics give
\begin{equation}\label{fbound}
 \sup_{t\le T}\|f_k(t)\|_\infty\le\|f_0\|_\infty\exp\!\left(\int_0^T \rho_*(t)\dd t\right) .
\end{equation}
With $G_k=\rho_k(u_k-v)f_k$ and $F_k=\int G_k\dd v$, mass and the second moment give
\begin{equation}\label{weighted}
 \int_0^\delta\iint\rho_k f_k\le C\int_0^\delta \rho_*(s)\dd s,\qquad
 \int_0^T\iint_{|v|>H}\rho_k f_k\le CH^{-2}\int_0^T \rho_*(s)\dd s .
\end{equation}
On bounded velocities away from zero time $\rho_k f_k$ is bounded. Thus \eqref{weighted} and
$\int_A|G_k|\le C(\int_A\rho_k f_k)^{1/2}$ give uniform integrability of $G_k,F_k$ and tight velocity tails for $G_k$. Write $e_k=\int v^2 f_k\dd v$. Splitting at velocity $H$ gives $n_k\le C_T H+H^{-2}e_k$ and $|j_k|\le C_T H^2+H^{-1}e_k$. Taking $H=e_k^{1/3}$ where $e_k>0$ gives $n_k^3+|j_k|^{3/2}\le C_T e_k$. Finally, \eqref{base} and \eqref{rhok} imply
\begin{equation}\label{pressuregain}
 \int_0^t\int(\rho_k^{\gamma+1}+\rho_k^2+|m_k|^2)\dd x\dd s
 \le C\int_0^t \rho_*(s)\dd s .
\end{equation}
The right side vanishes uniformly as $t\downarrow0$. On $[\delta,T]$ the density is bounded. Together with $u_{k,x}$ bounded in $L^2$, this gives uniform integrability of $p_k,\rho_k p_k,\rho_k u_{k,x}$ on the whole cylinder.

\subsection{Compactness and identification of the pressure}\label{densitycompact}
To pass to the limit in the convective flux, we need space-time compactness of the fluid momentum. The equations bound $\partial_t(\rho_k,m_k)$ in $L^1_t W^{-q,1}_x$ for $q$ large. The energy bounds these sequences in $L^\infty_tL^\gamma_x$ and $L^\infty_tL^{2\gamma/(\gamma+1)}_x$. The compact embedding of $L^p(\T_\ell)$ into $H^{-1}(\T_\ell)$ for $p>1$ and the time compactness theorem of \cite{Simon} give
\begin{equation}\label{compact}
 \rho_k\to\rho,\quad m_k\to m\quad\text{in }L^2_t H^{-1}_x,
 \qquad u_k\rightharpoonup u\quad\text{in }L^2_tH^1_x .
\end{equation}
Pairing with $u_k$ identifies $m=\rho u$ and $m_ku_k\to mu$ in distributions. On compact velocities $G_k=(m_k-\rho_kv)f_k$ is bounded in $L^2$ by \eqref{pressuregain} and \eqref{fbound}. For a smooth time-velocity cutoff $\psi$, $(\partial_t+v\partial_x)(\psi f_k)=-\partial_v(\psi G_k)+\psi_vG_k+\psi_tf_k$. All three coefficients are bounded in $L^2$, so velocity averaging \cite{DLM} gives compactness of the localized averages. Their $L^\infty$ bound removes the time-endpoint layers and gives strong convergence in every finite $L^p_{t,x}$. Pair with $m_k,\rho_k$ weak in $L^2$, then use \eqref{weighted}, to obtain
\[
 G_k\rightharpoonup\rho(u-v)f,\quad F_k\rightharpoonup F=\rho(nu-j)
     \quad\text{in }L^1 .
\]
Here $f_k\stackrel{*}{\rightharpoonup}f$ in $L^\infty$ on finite time intervals. The second moment bounds the mass and first-moment tails by $CH^{-2}$ and $CH^{-1}$, respectively, so $n_k\to n$ and $j_k\to j$ in $L^1_{t,x}$.

It remains to identify the weak limit of $p_k=\rho_k^\gamma$. The pressure can be paired with the density by testing momentum with the primitive of $\rho_k-r_k$. We use the following elementary estimate on an interval or a torus of length $\ell$. For $0<a\le\ell$, $0<\alpha\le1$ and $1\le p<\infty$, averaging $g\in C^{0,\alpha}$ over a subinterval of length $a$ containing $x$ gives
\begin{equation}\label{primitivecompact}
 \|g\|_\infty\le a^\alpha[g]_{C^{0,\alpha}}+a^{-1/p}\|g\|_p .
\end{equation}
In particular, the $L^\infty_tW^{1,\gamma}_x$ bound for $P_k=\I[\rho_k]$ gives a common spatial H\"older exponent $1-1/\gamma$, and $(P_k)_t=-m_k$ gives $\|P_k(t)-P_k(s)\|_2\le C|t-s|^{1/2}$ by \eqref{pressuregain}. Apply \eqref{primitivecompact} to this difference, first taking $a$ small, then $|t-s|$ small. Arzel\`a--Ascoli gives uniform convergence along a subsequence on the closed cylinder, including $t=0$. The spatial derivative and zero mean identify the limit as $P=\I[\rho]$. Test momentum by $\eta(t)P_k$, $\eta\in C_c^\infty(0,T)$. Since $m_kP_{k,t}+m_ku_kP_{k,x}=-r_km_ku_k$,
\begin{equation}\label{flux}
 \int\eta(p_k\rho_k-\nu\rho_k u_{k,x})
 =-\int\eta' m_kP_k+r_k\int\eta(m_ku_k+p_k)+\int\eta F_kP_k .
\end{equation}
All right-hand terms converge. In the limiting momentum equation with weak pressure limit $\bar p$, mollifications of $P$ converge uniformly, their time derivatives in $L^2$, and their spatial derivatives almost everywhere with a uniform bound on $\operatorname{supp}\eta$. Subtracting this test from \eqref{flux} yields the effective flux identity \cite{FNP}
\begin{equation}\label{cov}
 \int(\overline{\rho^{\gamma+1}}-\bar p\rho)
   =\nu\int(\overline{\rho u_x}-\rho u_x)\ge0\qquad\text{a.e. }t>0 .
\end{equation}
The density Young measure gives the left side as the integral of the covariance of $q^\gamma$ and $q$, nonnegative by monotonicity. Estimate \eqref{pressuregain} permits integration down to zero.

The sign in \eqref{cov} is what rules out oscillation of the density. To see this, use the renormalization $\rho\log\rho$ in the continuity equation. Since $\rho,u_x\in L^2_{t,x}$, spatial convolution in continuity has an $L^1$ commutator \cite{DL} tending to zero and retains the convolved initial value. Take $\beta(0)=0$ and $\beta'(q)=1+\log(\varepsilon+\min(q,R))$. Then
\[
 0\le q\beta'(q)-\beta(q)=\int_0^{\min(q,R)}\frac{z}{z+\varepsilon}\dd z\le q .
\]
As $R\to\infty$ and $\varepsilon\downarrow0$, the initial entropies converge to $\int\rho_0\log\rho_0$ by $\rho_0\in L^\gamma$, and the source converges to $\rho u_x$ in $L^1$. Subtracting the limiting entropy identity from its approximate counterpart, with initial convergence by \eqref{approx}, gives
\[
 0\le\int(\overline{\rho\log\rho}-\rho\log\rho)(t)
 =-\int_0^t\int(\overline{\rho u_x}-\rho u_x)\le0
 \qquad\text{a.e. }t .
\]
Testing the integrated balances up to zero retains the common initial entropy because \eqref{pressuregain} gives $\int_0^\delta\!\int|\rho_k u_{k,x}|\le C(\delta+\delta^{1-1/\gamma})^{1/2}$. Strict convexity and uniform integrability yield $\rho_k\to\rho$ a.e. along a subsequence and in $L^\gamma_{t,x}$, hence $p_k\to\rho^\gamma$ in $L^1$. Strong convergence of $\sqrt{\rho_k}$ in $L^2$ identifies the weak limit of $\sqrt{\rho_k}u_k$. Pairing $m_k,u_k$ in \eqref{compact} gives convergence of its squared norms, hence $\sqrt{\rho_k}u_k\to\sqrt\rho u$ in $L^2_{t,x}$.

\subsection{Particle labels and strong kinetic convergence}\label{kineticconstruct}
Lemma~\ref{masslemma} with $M_k=r_kx+P_k$ bounds the velocities uniformly even when $r_k\downarrow0$. Pushing forward $f_{0,k}$ gives, for $H$ sufficiently large,
\begin{equation}\label{tails}
 \sup_{t\ge0}\iint_{|v|>H}v^2 f_k(t)\dd x\dd v
 \le4\iint_{|v|>H/2}v^2 f_{0,k}\dd x\dd v .
\end{equation}
The right side tends to zero uniformly by \eqref{approx}. For the limiting equation $M=rx+P$ is continuous, $M_t=-m,M_x=\rho$ are locally in $L^2$ up to time zero, and $\|\rho\|_\infty\in L^1(0,T)$. The lemma gives \eqref{kinren} with trace $\beta(f_0)$ for each admissible renormalization $\beta$, and propagates finite initial velocity moments by \eqref{position}.

To obtain strong convergence of $f_k$, let $g$ be a weak-* limit of $f_k^2$. Then $f^2\le g\le C_Tf$. In $(\iint f_k^2)'=\iint\rho_k f_k^2$, strong convergence of $\rho_k$ permits passage to the limit on bounded velocities. Both tails are $O_T(H^{-2})$ by \eqref{fbound}, \eqref{weighted} and the second moment. The initial squares converge in $L^1$. Integrating \eqref{kinren} for $f^2$ is legitimate since $\rho|u-v|f^2\in L^1$. Subtraction gives an absolutely continuous defect $d=\iint(g-f^2)$ with

\[
 d(0)=0,\qquad d'(t)=\iint\rho(g-f^2)\le\|\rho(t)\|_\infty d(t).
\]
Gronwall gives $g=f^2$, hence strong $L^2$ convergence on bounded velocities. The moment tails \eqref{tails} give \eqref{strongf}.

We next identify the initial traces and the energy inequality for the limiting solution. The strong convergences yield $E_k\to E$ in $L^1(0,T)$ and along a subsequence at every time in a full-measure set $G$. Adjoin $0$ by \eqref{approx}. Viscosity is weakly lower semicontinuous. For the drag use
\[
 \int\zeta\rho_k f_k|u_k-v|^2
 \ge\int\zeta\rho_k f_k\{2(u_k-v)\xi-\xi^2\}
\]
with smooth $\zeta\ge0$ and smooth compactly velocity-supported $\xi$. The force and density convergences pass the right side to the limit. Taking the supremum over $\xi$ and approximating $\mathbf1_{(s,t)}$ from below by time cutoffs in \eqref{base} give
\begin{equation}\label{strongenergy}
 E(t)+\int_s^t D(\tau)\dd\tau\le E(s),\qquad s,t\in G,\quad s\le t .
\end{equation}
Pairings against smooth tests are absolutely continuous by the equations. The uniform fluid bounds, kinetic $L^\infty$ bound and tight mass tails give the weakly continuous representatives with the prescribed initial data. Lower semicontinuity extends \eqref{strongenergy} to every terminal time, and weak continuity extends \eqref{densitybound} to every $t>0$. The particle flow preserves mass.

Oddness gives zero trace for $u$. Specular tests extend across the endpoints by $(x,v)\mapsto(-x,-v)$ to continuous, piecewise $C^1$ tests. Periodic mollification keeps derivatives bounded and gives \eqref{spec} and \eqref{kinren} by dominated convergence. Even continuity tests and odd extensions of zero-boundary momentum tests give the fluid equations on $I$. Divide \eqref{strongenergy} by two and diagonalize in $T$. When $r=0$, mass conservation gives $\rho=m=0$, hence $\nu u_{xx}=0$ and $u=0$. The same approximations then converge in \eqref{strongf} to free specular transport.

\subsection{The energy equality at every time}\label{sec-energy-equality}
For the uncoupled compressible fluid, energy criteria have been proved by Yu \cite{Yu} and, with physical boundaries, Chen, Liang, Wang and Xu \cite{CLWX}. The following one-dimensional chain rule applies on the periodic lift.

\begin{lemma}[Kinetic energy at vacuum]\label{weightedchain}
Let $V=H^1(\T_\ell)$ and let $J$ be a bounded positive-time interval. Suppose $\rho\ge0$ is bounded, $m=\rho u\in L^\infty(J;L^2)$ and $u\in L^2(J;V)$. If $\rho_t\in L^\infty(J;V^*)$ and $m_t\in L^2(J;V^*)$, then the continuous $V^*$ representatives satisfy
\begin{equation}\label{weightedchainrule}
 \left(\int\frac{|m|^2}{2\rho}\right)'=
 \langle m_t,u\rangle_{V^*,V}-\tfrac12\langle\rho_t,u^2\rangle_{V^*,V}.
\end{equation}
The energy on the left is locally absolutely continuous at every interior time.
\end{lemma}
\begin{proof}
We take time differences of $\rho$ and $m$, whose derivatives are controlled by the equations. Write $\Delta_h g=(g(t+h)-g(t))/h$ and $u_h=u(t+h)$. Polarizing the kinetic energy at these two times gives, for a.e. $t$,
\[
 \Delta_h\left(\tfrac12\int\rho u^2\right)
 =\tfrac12\langle\Delta_h m,u_h+u\rangle-\tfrac12\langle\Delta_h\rho,u_hu\rangle .
\]
Since $V$ is an algebra, $u_hu\to u^2$ in $L^1_{\rm loc}(J;V)$. The differences $\Delta_h\rho$ are bounded in $L^\infty_{\rm loc}(J;V^*)$ and converge almost everywhere to $\rho_t$ in $V^*$, while $\Delta_h m\to m_t$ in $L^2_{\rm loc}(J;V^*)$. Testing in time gives \eqref{weightedchainrule} for an absolutely continuous representative $k$.

To identify its value at every $t$, observe that
\[
 \int\frac{|m(t)|^2}{2\rho(t)}
 =\sup_{\phi\in V}\{\langle m(t),\phi\rangle-\tfrac12\langle\rho(t),\phi^2\rangle\}.
\]
The spatial bounds persist at every time by weak compactness. Completing the square and approximating $m/(\rho+\epsilon)$ in $L^2$ gives the formula as $\epsilon\downarrow0$. Lower semicontinuity gives the upper bound $k(t)$ on this supremum. For the reverse bound, average $\|u\|_V^2$ over $(t+h/2,t+h)$ and choose Lebesgue times $t_j\downarrow t$ such that $h_j\|u(t_j)\|_V^2\to0$, where $h_j=t_j-t$. Testing the supremum at $t$ with $\phi=u(t_j)$ yields a value which differs from $k(t_j)$ by at most
\[
 \sqrt{h_j}\|m_t\|_{L^2(t,t_j;V^*)}\|u(t_j)\|_V+
 C h_j\|\rho_t\|_{L^\infty(J;V^*)}\|u(t_j)\|_V^2 .
\]
Both terms tend to zero. Continuity of $k$ gives the desired lower bound.

\end{proof}

\begin{proof}[Proof of the all-time balance \eqref{allenergy}]
Fix $0<\tau<T$. By \eqref{densitybound} and the energy bound, $m=\rho u\in L^\infty(\tau,T;L^2)$ and $mu\in L^2(\tau,T;L^2)$. With $F=\int\rho(u-v)f\dd v$, Cauchy--Schwarz gives
\[
 \|F(t)\|_1^2\le \|\rho(t)\|_\infty M_f\iint\rho f|u-v|^2 .
\]
Thus $F\in L^2(\tau,T;L^1)\subset L^2(\tau,T;V^*)$ and $p=\rho^\gamma$ is bounded there. The equations imply $\rho_t=-m_x\in L^\infty(\tau,T;V^*)$ and $m_t\in L^2(\tau,T;V^*)$. By Lemma~\ref{weightedchain}, the fluid kinetic energy satisfies
\[
 \frac{\dd}{\dd t}\int\frac{|m|^2}{2\rho}
 =\int p u_x-\nu\int |u_x|^2-\int Fu .
\]
Indeed $\frac12\langle\rho_t,u^2\rangle=\int mu\,u_x$ cancels precisely the convective term in $\langle m_t,u\rangle$. For the internal energy, renormalization gives $(\int\rho^\gamma/(\gamma-1))'=-\int p u_x$ on this interval. To identify both endpoints, convolve the continuity equation in space, using the weakly continuous density. The spatial commutator tends to zero in $L^1$ since $\rho$ is bounded and $u\in L^2_tH^1_x$. The derivative of the renormalization is bounded on the range, and the mollified pressures converge in $L^2_{t,x}$. At each fixed endpoint the spatial mollifications converge in $L^\gamma$ to the chosen density representative.

For the particles use an even smooth velocity truncation of $v^2/2$ with derivative bounded by $C|v|$. It is an admissible specular test. On the entire interval $(0,T)$,
\[
 \int_0^T\iint \rho v^2 f
 \le \sup_{t\le T}\iint v^2 f(t)\int_0^T\|\rho(t)\|_\infty\dd t<\infty .
\]
The drag dissipation then gives $\rho f\,|v|\,|u-v|\in L^1$. Removing the truncation yields
\[
 \frac{\dd}{\dd t}\frac12\iint v^2 f=\iint \rho f (u-v)v .
\]
The characteristic representation and \eqref{position} give continuity of the second moment at every time by dominated convergence. In the general periodic case one uses the velocity bound following \eqref{zper}. Adding the three identities gives \eqref{allenergy} for all $0<s\le t$. Finally \eqref{strongenergy} with $s=0$ bounds the right limit of the total energy by $E(0)$, and weak lower semicontinuity at the prescribed initial data gives the reverse bound. Since $D$ is integrable at zero, passage to $s=0$ proves the assertion.
\end{proof}

\subsection{The corrected momentum at the weak limit}\label{weakprimitive}
On $[0,T]\times[0,L]$, the moment bounds give $\|n_k(t)\|_3+\|j_k(t)\|_{3/2}\le C_T$. Since $(N_k)_x=n_k$ and $(N_k)_t=-j_k$ at the approximation level,
\[
 |N_k(t,x)-N_k(t,y)|\le C_T|x-y|^{2/3},\qquad
 \|N_k(t)-N_k(s)\|_{3/2}\le C_T|t-s|.
\]
Apply \eqref{primitivecompact} on $[0,L]$ to $N_k(t)-N_k(s)$ with $p=3/2$. Choosing $a=|t-s|^{3/4}\le L$ gives $\|N_k(t)-N_k(s)\|_\infty\le C_T|t-s|^{1/2}$ for small time increments. Thus the primitives are compact in the uniform norm on the closed strip. Since $n_k\to n$ in $L^1_{t,x}$ and $N_k(t,0)=0$, their only limit is $N(t,x)=\int_0^x n(t,y)\dd y$. The same integral at $t=0$ follows from \eqref{approx}. In particular $N_t=-j$ holds without appealing to a boundary trace of $j$.

Now $q_k=m_k+\rho_kN_k\to q$ in distributions. For $u_kq_k=m_ku_k+m_kN_k$, the first term converges by \eqref{compact} and the second by uniform convergence of $N_k$ and weak $L^2$ convergence of $m_k$. Also $p_k\to p$ in $L^1$. By the energy and \eqref{pressuregain}, $q_k$ is bounded in $L^\infty_tL^{2\gamma/(\gamma+1)}_x\cap L^2_{t,x}$. Hence $b_k$ has a common spatial H\"older bound with exponent $\alpha=(\gamma-1)/(2\gamma)>0$. Its equation is
\begin{equation}\label{Bweak}
 (b_k)_t=\nu u_{k,x}-u_kq_k-p_k+A_k(t).
\end{equation}
Here $u_{k,x}\in L^2_{t,x}$ uniformly and the other terms are bounded in $L^\infty_tL^1_x$. Integration gives $\|b_k(t)-b_k(s)\|_1\le C_T(|t-s|^{1/2}+|t-s|)$. Estimate \eqref{primitivecompact} with $p=1$ gives uniform equicontinuity in time. A subsequence converges uniformly to the zero-mean primitive $b$ of $q$. Its initial value is the primitive of $m_0+\rho_0N_0$, since these approximate initial momenta converge in $L^1(I)$.

It remains to identify the spatial average $A_k$. The integral of $m_ku_k$ over $I$ is half its integral over the reflected torus and converges against time tests by \eqref{compact}. The integrals of $m_kN_k$ and $p_k$ pass to the limit by the convergences above. Equation \eqref{Bweak} therefore yields \eqref{scalarpressure}, with $b_t\in L^1$ and $b_x\in L^2$ up to time zero. The estimate $\|b(t)\|_\infty\le C\|q(t)\|_1$ is uniform in time.

\section{The spatial phase for general periodic data}\label{general-periodic}
On a torus without reflection symmetry the same estimates hold after subtracting the conserved mean velocity. There is also a spatial translation, which can be computed from the initial data.

\begin{theorem}[Periodic weak solutions]\label{periodic-thm}
Let \eqref{data} hold on $\T_\ell$ with $r=\av{\rho_0}>0$. Set $M_f=\iint f_0$ and
\[
 U_*=\frac{\int m_0+\iint vf_0}{\ell r+M_f},\qquad
 P_0=\I[\rho_0],\qquad S_0(a)=a+P_0(a)/r .
\]
There is a global periodic weak solution with the properties of Theorem~\ref{exist}, using $H^1(\T_\ell)$ for the fluid velocity. In the moving frame $(x-U_*t,v-U_*)$, with fluid velocity $u-U_*$, the energy estimate \eqref{exp}, the particle convergence \eqref{w2}, and all the marginal and support assertions of Proposition~\ref{relax} hold. The limiting particle position associated with $(a,v)$ is
\begin{equation}\label{phase-label}
 [S_0(a)+(v-U_*)/r+\theta_*]_\ell,\qquad
 \theta_*=-\frac{\iint (v-U_*+P_0(a))f_0(a,v)\dd a\dd v}{r(\ell r+M_f)} .
\end{equation}
The length of the vacuum set obeys \eqref{vaclim} in both frames. In the moving frame, \eqref{vacrate}--\eqref{vacbl} hold on $\T_\ell$. In both the comparison \eqref{vaclabel} and the normalized distribution \eqref{vacnormalized}, the mass coordinate is $[S_0+\theta_*]_\ell$. The corresponding weight is given at the end of the proof. In laboratory coordinates the comparison measure at time $t$ is the translate by $U_*t$ of this moving-frame limit.
\end{theorem}

\begin{proof}
Use moving coordinates throughout the proof. Approximate as in \eqref{approx} with $r_k\ge r/2$. The strong approximants conserve total momentum. For each of them the fluid mass average $u_\rho=\av{m}/r$ (using its own $r,m$) satisfies
\begin{equation}\label{weighted-poincare}
 |u_\rho|^2\le (\ell r)^{-1}\int\rho u^2,\qquad
 \|u-u_\rho\|_\infty\le C_\ell\|u_x\|_2 .
\end{equation}
Indeed $u(x)-u_\rho=(\ell r)^{-1}\int\rho(y)(u(x)-u(y))\dd y$. This supplies the finite-time $L^2_t H^1_x$ bound.

\smallskip\noindent\emph{Density and compactness.}
Without walls, choose a periodic primitive and fix its additive constant by requiring the corrected momentum to have zero mean. With $\bar n=\av{n_0}$ let
\begin{equation}\label{perprimitive}
 N_1=\I[n]-(\bar n/r)P,\quad c(t)=\frac{\av{m+\rho N_1}}{r},
 \quad N=N_1-c(t),\quad q=m+\rho N,\quad b=\I[q].
\end{equation}
First use these definitions for each strong approximation with its own masses. Since $P_t=-m+\av m$ and $\partial_t\I[n]=-j+\av j$, we have
\[
 (\partial_t+u\partial_x)N_1=un-j+\av j-(\bar n/r)\av m .
\]
Adding the balance for $\rho N_1$ to momentum and integrating in space gives $c'=\av j-(\bar n/r)\av m$. Hence $N_t+uN_x=un-j$, $\av q=0$, and \eqref{scalarpressure} holds with $A=\av{uq+p}$. The energy, the masses and $r_k\ge r/2$ bound $N,c,b,A$ uniformly. Equation \eqref{logODE} gives \eqref{rhok} and the initial-layer bounds \eqref{weighted} and \eqref{pressuregain}.

In the effective flux test \eqref{flux} there is one extra term, $-\ell\int\eta\,\av{m_k}^2$, since $P_{k,t}=-m_k+\av{m_k}$. It converges by \eqref{compact}, which gives $\av{m_k}\to\av{m}$ in $L^2_t$. The covariance and entropy arguments therefore give the same strong convergence of the densities and of $\sqrt{\rho_k}u_k$.

Put $\Theta(t)=\int_0^t\av{m}(s)\dd s$ and $M=rx+P-\Theta$. Then $M_t=-m$, $M_x=\rho$. Along a trajectory $w=v+M$ is constant, and $z=v+P$ satisfies
\begin{equation}\label{zper}
 \dot z+rz=rP(t,X)+\av{m}(t).
\end{equation}
The analogous approximate equation, $r_k\ge r/2$ and the energy bound on $\av{m_k}$ give $|V_k(t)|\le e^{-r_kt}|v_0|+C$. This yields \eqref{tails}. Lemma~\ref{masslemma} applies to $M$, and the defect argument of Appendix~\ref{kineticconstruct} gives \eqref{strongf}. The velocity bound and Appendix~\ref{sec-energy-equality} give \eqref{allenergy} in this frame. Testing momentum by a spatial constant and the kinetic equation by truncated $v$ gives conservation of total momentum in the limit, since the force is integrable and the first-moment tails vanish uniformly. Thus $\av{j}=-\av{m}$. The energy in the laboratory frame differs by a constant fixed by mass and total momentum, so it also satisfies \eqref{allenergy}.

\smallskip\noindent\emph{Decay and the limiting phase.}
In \eqref{J}, $-r\int\rho u^2$ becomes $-r\int\rho(u-u_\rho)^2$ by the formula for $P_t$. This is controlled by viscosity using \eqref{weighted-poincare}. The extra term in the identity for $K$ is
\[
 2\av{m}\iint zf=\frac{2}{\ell}\left(\iint Pf-\iint zf\right)\iint zf
 \le\ell^{-1}\left(\iint Pf\right)^2 .
\]
Thus $K'\le-r K+C\|\rho-r\|_2^2$. Also $\ell\av{m}^2/r\le C\iint v^2f$ by zero total momentum. The proof of \eqref{exp} now applies to $\mathcal L=\E-\delta J+\delta^2K$. To control the fluid kinetic energy, use \eqref{weighted-poincare} on $u-u_\rho$ and the preceding bound on the mean.

The first moment of \eqref{zper} gives $(\iint zf)'=-r\iint vf+M_f\av{m}=(\ell r+M_f)\av{m}$. The velocity bound above justifies differentiation under the initial measure. Hence
\begin{equation}\label{phase-id}
 (\ell r+M_f)\Theta(t)=\iint (v+P(t,x))f(t,x,v)\dd x\dd v
                         -\iint(v+P_0(x))f_0(x,v)\dd x\dd v .
\end{equation}
The first moment on the right has square bounded by $M_f K(t)$. In the original initial variables this gives $\Theta_\infty/r=\theta_*$ in \eqref{phase-label}, with $|\Theta(t)-\Theta_\infty|^2\le C\E(t)$. Set $d=\Theta(t)-\Theta_\infty$. Then $\iint zf=(\ell r+M_f)d$ and
\[
 \iint(z-d)^2 f= K-(2\ell r+M_f)d^2\le K .
\]
Since $X=(w+\Theta-z)/r$, this gives the same bound \eqref{w2} for the translated limiting position. For a fixed label $y=w/r$ the position solves $\dot X=ry+\Theta-rX-P(t,X)$. Its derivative in $y$ still satisfies \eqref{jac}. On bounded labels the positions converge uniformly and exponentially to $y+\theta_*$ by \eqref{zper}, and the proof of \eqref{jacclose} applies. Comparing with this translation yields the marginal, velocity-support and moment assertions of Proposition~\ref{relax}.

\smallskip\noindent\emph{The vacuum measure.}
In \eqref{perprimitive}, $\I[n_k]$ converges uniformly as in Appendix~\ref{weakprimitive} since its time derivative is $-j_k+\av{j_k}$. The means $\av{j_k},\av{m_k}$ converge in $L^2_t$ on finite intervals by strong $L^1$ convergence of $j_k$, the moment bound and \eqref{compact}. The equation for $c_k'$ and convergence of $c_k(0)$ thus give uniform convergence of $c_k$. Passing to the limit as in \eqref{Bweak} proves the same pressure identity for the periodic weak solution.

Since $\av{N_1}=0$, $r|c|\le|\av m|+\|N_1\|_\infty\av{|\rho-r|}$. This decays exponentially, as does $P$. Thus $b-r\I^2[n]$ tends to zero uniformly and exponentially. The periodic kernel of $\I^2$ is Lipschitz, since its second derivative is $\delta_0-1/\ell$. Testing \eqref{w2} against its translates gives the same convergence of $b$ to $b_\infty=r\I^2[n_\infty]$. For $M_f=0$ the kernel term vanishes. The transported coordinate $S=M/r$ tends to $x-\theta_*$ modulo $\ell$. Apply \eqref{vacidentity} on the torus with $\varphi(\cdot+\theta_*)$. This gives \eqref{vacbl} and \eqref{vacnormalized} with $S_0(a)$ replaced by $[S_0(a)+\theta_*]_\ell$ and $W_0(a)=\exp((b_\infty([S_0(a)+\theta_*]_\ell)-b(0,a))/\nu)$, where $b(0)$ is given by \eqref{perprimitive}. Translation preserves the length and the Lipschitz norm of a test, proving the laboratory-frame assertions.

\end{proof}

The phase in \eqref{phase-label} is the same for every solution supplied by this construction with the given data. Continuous dependence of the full solution at vacuum is a separate question. Hoff and Serre \cite{HS} showed that it can fail for the uncoupled one-dimensional cavity problem in their solution class.

\section*{Declaration of generative AI use}
The author used OpenAI's ChatGPT to help draft and revise parts of the exposition, locate relevant literature, explore proof arguments and check calculations. This included the derivation of the relative weights of the vacuum set in Section~\ref{vacuum}. A separate AI-assisted reading of the draft suggested expanding the weak tests and compactness arguments. Responsibility for the mathematical arguments and references rests with the author.

\end{document}